\documentclass[12pt]{amsart}
\usepackage{epsfig,amssymb,latexsym,amsmath,graphics}

\newcommand{\bbz}{{\mathbb{Z}}}
\newcommand{\bbr}{{\mathbb{R}}}

\newcommand{\bbq}{{\mathbb{Q}}}
\newcommand{\bbc}{{\mathbb{C}}}

\newtheorem{definition}{Definition}[section]
\newtheorem{theorem}[definition]{Theorem}
\newtheorem{lemma}[definition]{Lemma}
\newtheorem{conjecture}[definition]{Conjecture}

\title[Vanishing of L-functions over number fields]{Vanishing of L-functions of elliptic curves over number fields}
\author{Chantal David}%
\address{Department of Mathematics and Statistics,
Concordia University,
1455 de Maisonneuve Blvd. West, Montr\'{e}al (Qu\'{e}bec), Canada H3G 1M8}%
\email[David]{cdavid@mathstat.concordia.ca}%
\author{Jack Fearnley}%
\email[Fearnley]{jack@mathstat.concordia.ca}%
\author{Hershy Kisilevsky}%
\email[Kisilevsky]{kisilev@mathstat.concordia.ca}%
\thanks{The first and third authors are partially supported by grants
from NSERC and FCAR}%
\begin{document}

\begin{abstract}
Let $E$ be an elliptic curve over $\bbq$, with L-function
$L_E(s)$. For any primitive Dirichlet character $\chi$, let
$L_E(s, \chi)$ be the L-function of $E$ twisted by $\chi$. In this
paper, we use random matrix theory to study vanishing of the
twisted L-functions $L_E(s, \chi)$ at the central value $s=1$. In
particular, random matrix theory predicts that there are
infinitely many characters of order 3 and 5 such that $L_E(1,
\chi)=0$, but that for any fixed prime $k \geq 7$, there are only
finitely many character of order $k$ such that $L_E(1, \chi)$
vanishes. With the Birch and Swinnerton-Dyer Conjecture, those
conjectures can be restated to predict the number of cyclic
extensions $K/ \bbq$ of prime degree such that $E$ acquires new
rank over $K$.
\end{abstract}

\maketitle

\section{Introduction}

Let $E$ be an elliptic curve defined over $\bbq$ with conductor
$N_E$. For any number field $K/\bbq$, let $E(K)$ be the group of
points of $E$ defined over $K$. By the Mordell-Weil Theorem,
$E(K)$ is a finitely generated abelian group.
Let $L_E(s, K)$ be the L-function of $E$ over the field $K$.

\begin{conjecture}[Birch and Swinnerton-Dyer conjecture over number fields]
$L_E(s, K)$ has analytic continuation to the whole complex plane, and
\begin{eqnarray*}
\displaystyle \mbox{ord}_{s=1} L_E(s, K) &=& r_K(E)
\end{eqnarray*}
where $r_K(E)$ is the rank of $E(K)$.
\end{conjecture}

In this paper, we fix $E$ an elliptic curve over $\bbq$, and we
study how  the rank varies over abelian fields $K / \bbq$ of fixed prime
degree. For example, are there infinitely many such number fields
where $E$ acquires new rank over $K$ (i.e. $r_K(E) > r_\bbq(E)$)? With the Birch and
Swinnerton-Dyer conjecture, one can rephrase the question in terms
of vanishing of the L-function $L_E(s,K)$ at $s=1$. Let $K$ be an
abelian extension of $\bbq$ with Galois group $G$ and conductor
$m$. Let $\hat{G}$ be the group of characters of $G$ which can be identified with a set
of Dirichlet characters
$$
\chi: \left( \bbz / m \bbz \right)^* \rightarrow \bbc^* .
$$
Let
$$L_E(s) = L_E(s, \bbq) = \sum_{n \geq 1} \frac{a_n}{n^s}$$
be the L-function of $E$ over $\bbq$. For each primitive Dirichlet character
$\chi$, let
$$L_E(s, \chi) = \sum_{n \geq 1} \frac{\chi(n) \, a_n}{n^s}$$
be the L-function of $E$ over $\bbq$ twisted by the character
$\chi$. By the work of \cite{wiles, taylorwiles, BCDT}, $L_E(s)$ and
$L_E(s, \chi)$ have analytic continuation to the whole complex plane.
It also follows from properties of number fields that
\begin{eqnarray}
\label{numberfields}
 L_E(s, K) = \prod_{\chi \in \hat{G}} L_E(s, \chi) ,
\end{eqnarray}
and the vanishing of the twisted L-functions $L_E(s, \chi)$ at
$s=1$ is
equivalent, via the Birch and Swinnerton-Dyer conjecture, to the
existence of rational points of infinite order on $E(K)$.

In this paper, we use random matrix theory to
study the vanishing of the twisted L-functions
$L_E(s, \chi)$ at $s=1$. It has been known
since the work of Montgomery \cite{montgomery} that certain
statistics on probability spaces of random matrices (as pair
correlation between the eigenangles of the matrices) are similar
to the same statistics on the zeroes of the Riemann zeta function.
This intuition is supported by the extensive computations of
Odlyzko \cite{odlyzko} on the critical zeroes of the Riemann zeta
function.

This was explored further in the work of Katz and
Sarnak \cite{KaSa1, KaSa2}, who extend the analogy between other
probability spaces of matrices, and families of L-functions. Katz
and Sarnak also studied the case of function fields, where they
can actually prove some of those mysterious connections between
random matrices and families of L-functions.

In order to study different statistics of number theoretic
objects, Keating and Snaith \cite{KS1, KS2} introduced a new
random variable on spaces of random matrices, the characteristic
polynomial of the matrix evaluated at a given point. They computed
the probability distribution of this new variable, which led to
striking conjectures for the asymptotic behavior of the moments of
the Riemann zeta function on the critical line.
%%Before the work of
%%Keating and Snaith, even this conjectural asymptotic could not be
%%obtained by purely number theoretic methods.
The ideas of Keating
and Snaith have been applied to study vanishing of L-functions in
families \cite{CKRS, DFK, watkins}. Families of quadratic twists
are studied in \cite{CKRS}, where a conjectural asymptotic for
the number of quadratic twists with even non-zero rank is
presented. This refines a conjecture of Goldfeld \cite{goldfeld}
which predicts that quadratic twists with rank greater than one
have density zero. In \cite{DFK}, the authors used the ideas of
Keating and Snaith to obtain a conjectural asymptotic for the
number of cubic Dirichlet characters $\chi$ such that $L_E(1,
\chi)$ vanishes. We present in this paper the case of
characters of order $k$, for $k$ any odd prime. More precisely,
the conjectures that we obtain from the random matrix model are
(the case $k=3$ of \cite{DFK} is included for completeness)

\begin{conjecture} \label{mainconjecture}
Let $k$ be an odd prime, let $E$ be an elliptic curve defined over
$\bbq$, and let
\begin{eqnarray*}
N_{E,k}(X)  &=& \# \left\{ \chi \;\mbox{of order $k$} \;:\;
\mbox{cond}(\chi) \leq X \;\mbox{and}\; L_E(1 ,\chi) = 0
\right\} .
%%%N_{E,k}(X) &=& \# S_{E,k}(X)
\end{eqnarray*}

If $k = 3$, then
\begin{eqnarray*}
\log{N_{E,k}(X)} &\sim&  \frac{1}{2} \log{X} \;\;\;\;
\mbox{as $X \rightarrow \infty$.}
\end{eqnarray*}

If $k = 5$, then
$N_{E,k}(X)$ is unbounded, but
$N_{E,k}(X) \ll X^\epsilon$ for any $\epsilon > 0$
as $X \rightarrow \infty$.

If $k \geq 7$, then
$N_{E,k}(X)$ is bounded.

\end{conjecture}

In the light of (\ref{numberfields}), and
under the Birch and Swinnerton-Dyer conjecture,
one can rewrite $N_{E,k}(X)$ as
\begin{eqnarray*}
\label{twistsandfields}
N_{E,k}(X) &=& (k-1) \;  \# \left\{ K/\bbq \; \mbox{cyclic of degree $k$} \;:\; \right.\\
&& \;\;\;\; \left. \mbox{cond}(K) \leq X  \;\mbox{and}\; r_K(E)
> r_\bbq(E) \right\}.
\end{eqnarray*}

The structure of the paper is as follows. In the second section,
we use modular symbols to rewrite the special values $L_E(1,
\chi)$ as a product of terms depending only on $E$, and some
algebraic integer $n_E(\chi)$ depending on the character. In the
third section, we use the embedding of number fields as lattices
in $\bbc$ to give a discretisation of the algebraic integers
$n_E(\chi)$. In the fourth section, we use this discretisation and
the work of Keating and Snaith to obtain conjectures on the
asymptotic behavior of $N_{E,k}(X)$. Finally, the last section
presents some experimental results.

\section{Special values and modular symbols}

The notation of this section follows the introduction of
\cite{MTT}. Let $E$ be an elliptic curve over $\bbq$, and let
$f(z) = \sum_{n \geq 1} a_n \; e^{2 \pi i n z}$ be the Fourier
expansion of $E$. Then, the L-function $L_E(s)$ is the Mellin
transform
\begin{equation}
\label{mellintransform} L_E(s) = \frac{(2 \pi)^s}{\Gamma(s)}
\int_{0}^{\infty} f (it) t^{s-1} \; dt .
\end{equation}
By the work of \cite{wiles, taylorwiles, BCDT}, $f$ is modular and
$L_E(s)$ has analytic continuation to the whole complex plane, and
satisfies the functional equation
\begin{equation}
\label{functionalequation} \Lambda_E(s) = \left(
\frac{\sqrt{N_E}}{2 \pi} \right)^s \Gamma(s) L_E(s) = w_E
\Lambda_E(2-s)
\end{equation}
where $w_E = \pm 1$ is called the root number.
From (\ref{mellintransform}), we have that $L_E(1) = 2 \pi
\int_{0}^{\infty} f(it) \; dt.$ For $a,m \in \bbq, m> 0$, one
defines the modular symbols
\begin{eqnarray}
\label{defmodularsymbol} \lambda(a,m;E) &=& \lambda(a,m;f) = 2 \pi
\int_{0}^{\infty} f \left( it - \frac{a}{m} \right) dt.
\end{eqnarray}
Let $\chi$ be a primitive character of modulus $m$ with Gauss sum
$$\tau(\chi) = \sum_{a \;\rm{mod}\;  m} \chi(a) e^{2 \pi i a/m}. $$
The twisted L-function $L_E(s,\chi)$ satisfies the functional
equation
\begin{eqnarray*}
\label{generalFE} \Lambda_E(s, \chi) &=& \left( \frac{m
\sqrt{N_E}}{2\pi} \right)^s \Gamma(s) L_E(s, \chi) \\
&=&\frac{w_E \chi(N_E) \tau(\chi)^2}{m} \Lambda_E(2-s,
\overline{\chi}).
\end{eqnarray*}
From the identity
$$
\chi(n) = \frac{1}{\tau(\overline{\chi})}\sum_{a \;\rm{mod}\;  m}
\overline{\chi}(a) e^{2 \pi i a n /m},
$$
we have
\begin{eqnarray*}
f_\chi(z) &=& \sum_{n \geq 1} \chi(n) a_n e^{2 \pi i n z}\\
&=& \frac{1}{\tau(\overline{\chi})}\sum_{a \;\rm{mod}\;  m}
\overline{\chi}(a) f \left( z + a/m \right)
\end{eqnarray*}
by rearranging the sums (Birch's lemma).
%%It then follows that for $b,n \in \bbq, b> 0$, we have
%%\begin{eqnarray*}
%%\lambda(b, n; f_\chi) &=& \frac{1}{\tau(\overline{\chi})}\sum_{a
%%\mod m} \overline{\chi}(a) \lambda(mb-na, mn; f).
%%\end{eqnarray*}
%%Using $b=0, n=1$ and (\ref{defmodularsymbol}), we get
It then follows that
\begin{eqnarray}
\label{specialvalue1} L_E(1, \chi) &=&
\frac{1}{\tau(\overline{\chi})} \sum_{a \;\rm{mod}\;  m}
\overline{\chi}(a) \lambda(a,m;E).
\end{eqnarray}

We define
\begin{eqnarray*}
\lambda^+(a,m;E) &=& \lambda(a,m;E) + \lambda(-a,m;E) \\
&=& {2 \pi} \int_{0}^{\infty} \sum_{n \geq 1} a_n \;e^{-2 \pi n t}
\;\left( e^{2 \pi i a n/m} + e^{-2 \pi i a n / m} \right) dt \\
&=& 4 \pi \sum_{n \geq 1} a_n \left( \mbox{Re}\left\{ e^{2 \pi i a
n /
m} \right\} \right)\int_{0}^{\infty} e^{-2 \pi n t} dt \\
&=& 2 \sum_{n \geq 1} \frac{a_n}{n} \;\mbox{Re}\left\{ e^{2 \pi i
a n / m} \right\}
\end{eqnarray*}
which is a real number as $E$ is defined over $\bbq$. Then, all
$\lambda^+(a,m;E)$ can be written as a rational multiple of the
real period, which we denote by $\Omega_E$. We define
$$
\Lambda(a,m;E)  = \frac{\lambda^+(a,m;E)}{\Omega_E} \in \bbz.
$$
In this paper, $\chi$ has prime order $k \geq 3$. Then, $\chi(-1)
= 1$, and we can rewrite (\ref{specialvalue1}) as
\begin{eqnarray}
\label{specialvalue2} L_E(1, \chi) &=& \frac{\Omega_E}{2
\tau(\overline \chi)} \sum_{a \;\rm{mod}\;  m} \overline{\chi}(a)
\Lambda(a,m;E)
\end{eqnarray}
where the integers $\Lambda(a,m;E)$ do not depend on the character $\chi$,
but only on the conductor $m$.
Then, we see from (\ref{specialvalue2}) that
$$
\frac{2 \tau(\overline \chi) L_E(1, \chi)}{\Omega_E}
$$
is an algebraic integer in the field obtained by adding a $k$th root of unity.
%%\in \bbz[\xi_k].$$
In fact, one can prove the stronger result which is critical to
the  discretisation of the next section. The particular case $k=3$
was proven in \cite{DFK}. For any odd prime $k$,
let $\bbq(\xi_k)$ be the cyclotomic
field obtained by adding a primitive $k$th root of unity $\xi_k$,
and let $\bbq(\xi_k)^+$ be the maximal real extension $\bbq
\subseteq \bbq(\xi_k)^+ \subseteq \bbq(\xi_k)$ of degree
$(k-1)/2$ over $\bbq$. The ring of integers of $\bbq(\xi_k)^+$ will be
denoted by $\bbz[\xi_k]^+$.

\begin{theorem}
\label{maximalreal} Let $k$ be an odd prime, and let $\chi$ be a primitive character of order
$k$. Then,
\begin{eqnarray*}
\frac{2 \tau(\overline \chi) L_E(1, \chi)}{\Omega_E} = \left\{
\begin{array}{ll} \chi(N_E)^{(k+1)/2} \;n_E(\chi) &
 \mbox{when $w_E = 1$} \\
 \\
\displaystyle \left( \xi_k^{-1} - \xi_k \right)^{-1}
 \chi(N_E)^{(k+1)/2} \;n_E(\chi) &
 \mbox{when $w_E = -1$} \end{array} \right.
\end{eqnarray*}
where
$n_E(\chi) \in \bbz[\xi_k] \cap \bbr = \bbz[\xi_k]^+$.
\end{theorem}

\noindent{\bf Proof:} From the functional equation, we have
\begin{eqnarray*}
L_E^{alg}(1, \chi) &=& \frac{2 \tau(\overline \chi) L_E(1, \chi)}{\Omega_E} =
\frac{2 \tau(\overline \chi) w_E \chi(N_E) \tau(\chi)^2}{m \Omega_E} L_E(1, \overline \chi) \\
&=& w_E \chi(N_E) \overline{\frac{2 \tau(\overline \chi) L_E(1, \chi)}{\Omega_E}}
= w_E \chi(N_E) \overline{L_E^{alg}(1,\chi)}
\end{eqnarray*}
Let $z \in \bbc^*$ satisfying $z = w_E \chi(N_E) \overline{z}$.
Then, $L_E^{alg}(1,\chi) = \alpha \overline{z}^{-1}$ with $\alpha$ real.
If $w_E  = 1$, we take $z = \chi(N_E)^{(k+1)/2}$, and $L_E^{alg}(1,\chi) =
\alpha \overline{z}^{-1}$ with $\alpha \in \bbr \cap \bbz[\xi_k] = \bbz[\xi_k]^+$,
which gives the result. \\
If $w_E = -1$, we take $z = \left( \xi_k - \xi_k^{-1} \right)
\chi(N_E)^{(k+1)/2}$, and
$L_E^{alg}(1, \chi) = \alpha \overline{z}^{-1}$ with $\alpha \in \bbr \cap \bbz[\xi_k] = \bbz[\xi_k]^+$,
which gives the result.
\qed

\section{Discretisation}
\label{discretisation}

By Theorem \ref{maximalreal},
$n_E(\chi)$ is an algebraic integer in $\bbz[\xi_k]^+$, and there is then a natural
discretisation on the algebraic integer $n_E(\chi)$ given by the geometry of numbers.
Let $\phi$ be the map
\begin{eqnarray*}
\phi :   \bbz[\xi_k]^+ &\rightarrow& \bbr^{(k-1)/2} \\
\alpha &\mapsto& (\sigma_1(\alpha), \sigma_2(\alpha), \dots, \sigma_{(k-1)/2}(\alpha) ) \\
\end{eqnarray*}
where $\mbox{Gal}(\bbq(\xi_k)^+/\bbq) = \left\{ \sigma_1 = 1, \sigma_2, \dots, \sigma_{(k-1)/2} \right\}$.
Let $\alpha_1, \dots, \alpha_{(k-1)/2}$ be an integral basis for $ \bbz[\xi_k]^+$.
The image of $ \bbz[\xi_k]^+$ in $\bbr^{(k-1)/2}$ is the lattice generated by
the linearly independent vectors
$$\omega_1 = \phi(\alpha_1), \;\dots, \;\omega_{(k-1)/2} = \phi(\alpha_{(k-1)/2}).$$
Let $R \subseteq \bbr^{(k-1)/2}$ be the region
\begin{eqnarray*}
R &=& \left\{\; a_1 \omega_1 + a_2 \omega_2 + \dots + a_{(k-1)/2}
\omega_{(k-1)/2}: \right. \\
&& \;\;\; \left. -1 < a_i < 1 \;\mbox{for $1 \leq i \leq (k-1)/2$}
\;\right\}.
\end{eqnarray*}
The discretisation given by the embedding of
$\bbz[\xi_k]^+$ in $\bbr^{(k-1)/2}$ is then
\begin{eqnarray}
\label{geometryofnumbers} n_E(\chi) = 0 \iff \phi(n_E(\chi)) \in R.
\end{eqnarray}
Let $\chi$ be any character of conductor $m$ and order $k$. For
any automorphism $\sigma \in \mbox{Gal}(\bbq(\xi_k) / \bbq)$, let
$\chi^{\sigma}$ be the character
\begin{eqnarray*}
\chi^{\sigma} : \left( \bbz / m \bbz \right)^* &\rightarrow&
\left< \xi_k \right>
\subseteq \bbc^* \\
a &\mapsto& \sigma \left( \chi(a) \right)
\end{eqnarray*}
Then, $\chi^\sigma$ is also a character of conductor $m$ and order
$k$.

\begin{lemma}
\label{discrete}
Let $k$ be an odd prime, and $\chi$ a character of
order $k$ and conductor $m$. For any $\sigma$ in
Gal$(\bbq(\xi_k)/\bbq)$, we have
\begin{eqnarray*}
\vert L_E(1, \chi^\sigma) \vert &=& \frac{c_{E,k}}{m^{1/2}} \vert
n_E(\chi)^\sigma \vert
\end{eqnarray*}
where $c_{E,k}$ is an explicit constant depending only on $E$ and
$k$.
\end{lemma}
\noindent{\bf Proof:}  Using (\ref{specialvalue2}), we have
\begin{eqnarray*}
L_E^{alg}(1,\chi)^\sigma &=& \left( \frac{2 \tau(\overline \chi) L_E(1,
\chi)}{\Omega_E} \right)^\sigma \\
&=& \sum_{a \;\rm{mod}\; m} \overline{\chi}^\sigma(a)
\Lambda(a,m;E) = L_E^{alg}(1,\chi^\sigma).
\end{eqnarray*}
Suppose first that $\omega_E=1$. Then,
\begin{eqnarray*}
n_E(\chi)^\sigma &=& L_E^{alg}(1,\chi)^\sigma \left( \chi(N_E)^{(k+1)/2}
\right)^\sigma \\
&=& \frac{ 2 \tau(\overline{\chi}^\sigma) L_E(1,
\chi^\sigma)}{\Omega_E} \left( \chi(N_E)^{(k+1)/2} \right)^\sigma,
\end{eqnarray*}
and taking absolute values we get the result with $c_{E,k} =
\Omega_E / 2$. The proof for $\omega_E = -1$ is similar, with a
different explicit constant $c_{E,k}$. \qed

We first consider the case $k=5$. We have that
$\bbz[\xi_5]^+ = \bbz\left[ \alpha \right]$ with $\alpha = (1 +
\sqrt{5})/2$, and $G_5 = \left< 1, \tau \right>$, where the
non-trivial automorphism $\tau$ sends $\sqrt{5}$ to $-
\sqrt{5}$. Then, the lattice of $\bbz\left[ \alpha \right]$ in
$\bbr^2$ is generated by $\omega_1 = (\alpha, \alpha^\tau)$ and
$\omega_2 = (\alpha^\tau, \alpha)$.
Let $R$ be the region
$$
R = \left\{ a \omega_1 + b \omega_2 : -1 < a < 1, -1 < b < 1 \right\}.
$$
By (\ref{geometryofnumbers}), $n_E(\chi) = 0$ if and only if
$\phi(n_E(\chi)) = (n_E(\chi), n_E(\chi)^\tau) \in R.$ As the region
$R$ is not symmetric with respect to the absolute value, and we
have a probability model for $\left \vert L_E(1, \chi) \right
\vert$, we also consider the two regions of $\bbr^2$
\begin{eqnarray*}
R_1 &=& \left\{ (x,y) \; : \;  - 1 < x,y <   1 \right\} \\
R_2 &=& \left\{ (x,y) \; : \;  -\sqrt{5} < x,y <  \sqrt{5}  \right\}
\end{eqnarray*}
with the property that $R_1 \subseteq R \subseteq R_2$, and
\begin{eqnarray}
\label{absolutevalue}
 (n_E(\chi), n_E(\chi)^\tau) \in R_i &\iff&
\left( \vert n_E(\chi) \vert, \vert n_E(\chi)^\tau \vert \right)
\in \vert R_i \vert
\end{eqnarray}
where
$$
\vert R_i \vert = \left\{ (x,y) \in R_i \; : \; x,y \geq 0 \right\} .
$$
The following lemma is now immediate from  (\ref{absolutevalue})
and Lemma \ref{discrete}
\begin{lemma}
\label{fork5} Let $\sigma \in \mbox{Gal}(\bbq(\xi_5)/\bbq)$ be an
automorphism which restricts to the non-trivial automorphism of
$\bbq(\sqrt{5})$. For $i=1,2$, we have
\begin{eqnarray*}
\left( n_E(\chi), n_E(\chi)^\sigma \right) \in R_i &\iff&   \vert L_E(1, \chi)
\vert, \vert L_E(1, \chi^\sigma) \vert \leq \frac{c_i}{\sqrt{m}}
\end{eqnarray*}
where $c_1, c_2$ are explicit constants depending only on $E$.
\end{lemma}

We now suppose that $k \geq 7$. Let $B$ be the non-zero constant
$$
B = \max_{1 \leq i \leq (k-1)/2} \sum_{j=1}^{(k-1)/2} \vert \sigma_i ( \alpha_j)
\vert ,
$$
and let $R  \subseteq R' \subseteq \bbr^{(k-1)/2}$ be the region
$$
R' = \left\{ (x_1, \dots, x_{(k-1)/2} ) \; : \; -B \leq x_i \leq B
\;\mbox{for $1 \leq i \leq (k-1)/2.$} \right\} $$ Then, $n_E(\chi)
= 0  \Rightarrow \phi(n_E(\chi)) \in R'$ by
(\ref{geometryofnumbers}). The following lemma is now immediate
from  Lemma \ref{discrete}
\begin{lemma}
\label{fork7} Let $\sigma_1, \dots, \sigma_{(k-1)/2} \in
\mbox{Gal}(\bbq(\xi_k)/\bbq)$ be a set of representatives for
$G_k = \mbox{Gal}(\bbq(\xi_k)^+/\bbq)$. Then,
$\phi(n_E(\chi)) \in R'$ if and only if
\begin{eqnarray*}
\vert L_E(1, \chi^{\sigma_i}) \vert \leq \frac{c_k}{m^{1/2}} \;\;\;
\mbox{for $1 \leq i \leq (k-1)/2$}
\end{eqnarray*}
where $c_k$ is an explicit constant depending only $E$ and $k$.
\end{lemma}

\section{Unitary Random Matrices}

Let $U(N)$ be the set of unitary $N \times N$ matrices with complex
coefficients which forms a probability space with respect to the
Haar measure. For each $A \in U(N)$, let
\begin{eqnarray*}
P_A(\lambda) = \mbox{det} (A - \lambda I)
\end{eqnarray*}
be the characteristic polynomial of $A$. For any $s \in \bbc$, let
$$
M_N(s) = \int_{U(N)} |P_A(1)|^s \;d{\mbox {Haar}} $$ be the
moments for the distribution of $|P_A(1)|$ in $U(N)$ with respect
to the Haar measure. In \cite{KS1}, Keating and Snaith proved that
\begin{eqnarray}
\label{MUSN} M_U(s,N)  &=& \prod_{j=1}^{N} \frac{\Gamma(j)
\Gamma(j+s)}{\Gamma^2(j+s/2)} ,
\end{eqnarray}
and then $M_U(s,N)$ is analytic for $\mbox{Re}(s) > -1$, and has
meromorphic continuation to the whole complex plane. By Fourier inversion, the
probability density of $|P_A(1)|$ is
$$
p(x) = \frac{1}{2\pi i} \int_{(c)} M_N(s) x^{-s-1} ds
$$
for some $c > -1$. Then, for any $I  \subseteq \bbr$,
$$
\mbox{Prob} \left( |P_A(1)| \in I \right) = \int_{I} p(x) \; dx .
$$
In our application to the vanishing of twisted L-functions, we will
be interested only in small values of $x$ where  the value of
$p(x)$ is determined by the first pole of $M_U(s,N)$ at $s=-1$.
More precisely, for
$$x \leq N^{-1/2},$$
one can show that
$$
p(x) \sim G^{2}(1/2) N^{1/4} \;\;\;\;\; \mbox{as $N \rightarrow
\infty$},
$$
where $G(z)$ is the Barnes G-function, with special value
$$G(1/2) = \exp{\left( \frac{3}{2} \zeta'(-1) - \frac{1}{4} \log{\pi} + \frac{1}{24} \log{2} \right)}$$
(see \cite[p. 81]{KS1} or \cite[p. 58]{hughesthesis} for more details).

We now consider the moments for the special values of L-functions
in families of twists. Fix $k \geq 3$, and let
\begin{eqnarray*}
S_k(X) &=& \left\{ \chi \; \mbox{of order $k$ and conductor
$\leq X$} \right\} \\
N_k(X) &=& \# S_k(X) \sim b_k X
\end{eqnarray*}
with an explicit constant $b_k$ (see for example \cite{cohenetal}). We then define
for any $s \in \bbc$
\begin{eqnarray}
\label{allmoments} M_E(s,X) &=& \frac{1}{N_k(X)} \sum_{\chi \in
S_k(X)} |L_E(1, \chi)|^s
\end{eqnarray}
The family of twists of order $k$ has unitary symmetry, as the
values  $|\zeta(1/2 + it)|$ on the critical line. Then

\begin{conjecture} [Keating and Snaith Conjecture for twists of order $k$]
\label{conjectureKS}
\begin{eqnarray*}
M_E(s,X) &\sim& a_E(s/2) M_U(s,N)  \;\;\;\;\; \mbox{as $N = 2 \log{X} \rightarrow \infty$},
\end{eqnarray*}
where $a_E(s/2)$ is an arithmetic factor
depending only on the curve $E$.
\end{conjecture}
In the conjecture, the relation between $N$ and $X$ is obtained by
equating the mean density of eigenangles of matrices in the
unitary group, and the mean density of non-trivial zeroes of the
twisted L-functions $L_E(s, \chi)$ at a fixed height (see
\cite{DFK}). The arithmetic factor $a_E(s)$ can not be obtained
from the random matrix theory, and has to be determined separately
for each family from its arithmetic. This was done for the family
of cubic twists in \cite{DFK}, and could be done for the family of
twists of order $k$ for each $k$. The arithmetic factor $a_E(s)$
would then be a meromorphic function for all $s \in \bbc$. As it
will be seen below, the only influence of the arithmetic factor
$a_E(s)$ in our application is that the special value $a_E(-1/2)$
will be part of the constant of the conjectural asymptotic of
$N_{E,k}(X)$. This would not provide any further information to
the cases $k \geq 5$ considered in this paper considered in this
paper in view of Conjecture \ref{mainconjecture}.

>From Conjecture \ref{conjectureKS}, the probability density $p_E(x)$ for
the distribution of the special values $|L_E(1, \chi)|$ for
characters of order $k$ is
\begin{eqnarray}
\nonumber
p_E(x) &=& \frac{1}{2\pi i} \int_{(c)} M_E(s,X) \,x^{-s-1} \,ds \\
\label{computep} &\sim& \frac{1}{2\pi i} \int_{(c)} a_E(s/2)
\,M_U(s,N) \,x^{-s-1} \,ds
\end{eqnarray}
as $N = 2 \log{X} \rightarrow \infty$.
As above, when $x \leq N^{-1/2}$, the value of $p_E(x)$ is determined
by the residue of $M_U(s,N)$ at $s=-1$, and it follows from (\ref{computep})
that
\begin{eqnarray}
\label{smallx}
p_E(x) \sim C_E \log^{1/4}{X}
\end{eqnarray}
for $x \leq (2 \log{X})^{-1/2}$, $X \rightarrow \infty$, and
$C_E = 2^{1/4} a_E(-1/2) G^2(1/2)$.

Let $\chi$ be a character of order $k \geq 3$ and conductor $m$.
We apply the above model to find the probability that $|L_E(1, \chi)| < c m^{-1/2}$,
for some constant $c > 0$.
For $x < c m^{-1/2} < (2 \log{m})^{-1/2}$ (for $m$ large enough), we have
$p_E(x) \sim C_E \log^{1/4}{m}$, and then
\begin{eqnarray}
\nonumber
\mbox{Prob} \left( |L_E(1, \chi)| < c m^{-1/2} \right) &\sim&
\int_0^{c m^{-1/2}} C_E \log^{1/4}{m} \; dx \\
\label{probability} &=& c \,C_E \frac{\log^{1/4}{m}}{m^{1/2}} .
\end{eqnarray}

We now use the probability density of the random matrix model
with the discretisation of Section \ref{discretisation} to obtain
conjectures for the vanishing of the L-values $L_E(1, \chi)$. We
first suppose that $k=5$, and as in the previous section, let
$\sigma \in \mbox{Gal}(\bbq(\xi_5)/\bbq)$ which restricts to the
non-trivial automorphism of $\bbq(\sqrt{5})$. We saw in the
previous two sections that
\begin{eqnarray*}
L_E(1, \chi) = 0
&\iff& \phi(n_E(\chi)) = (n_E(\chi), n_E(\chi)^\sigma) \in R.
\end{eqnarray*}
As $R_1 \subseteq R \subseteq R_2$, and using Lemma \ref{fork5},
the probability that
$L_E(1, \chi)$ is zero is bounded below by
\begin{eqnarray*}
%%\label{quinticbelow}
\mbox{Prob} \left( \vert L_E(1, \chi) \vert <
\frac{c_1}{\sqrt{m}} \right) \mbox{Prob} \left( \vert L_E(1,
\chi^{\sigma}) \vert < \frac{c_1}{\sqrt{m}} \right)
\end{eqnarray*}
and bounded above by
\begin{eqnarray*}
%%\label{quinticabove}
\mbox{Prob} \left( \vert L_E(1, \chi) \vert <
\frac{c_2}{\sqrt{m}} \right) \mbox{Prob} \left( \vert L_E(1,
\chi^{\sigma}) \vert < \frac{c_2}{\sqrt{m}} \right) .
\end{eqnarray*}
Assuming that $|L_E(1, \chi)$ and $|L_E(1, \chi^\sigma)|$ are independent identically distributed
random variables, and using (\ref{probability}), we get that the probability that
$L_E(1, \chi)$ is zero is about
$$
\frac{\log^{1/2}{m}}{m} ,
$$
neglecting all constants which are not significant here. The sum
of the probabilities is
\begin{eqnarray}
\sum_{\chi \in S_3(X)} \frac{\log^{1/2}{m}}{m} &\sim& \frac{2
b_3}{3} \log^{3/2}{X} .
\end{eqnarray}
As discussed in \cite{DFK}, the exact power of $\log{X}$ that is
obtained with the random matrix approach depends subtly on the
discretisation, and is difficult to predict. For examples,
rational torsion of order three on the elliptic curve seemed to
cause extra vanishing of the twisted L-values $L_E(1, \chi)$ for
cubic characters, and changed the power of logarithm in the
conjectural asymptotic for $N_{E,k}(X)$ of \cite{DFK}. For $k =
5$, the sum of the probabilities is just on the border between
convergence and divergence, and the random matrix model seems to
indicate that the number of quintic twists such that $L_E(1,
\chi)$ vanishes is infinite, but that $N_{E,k}(X) \ll X^\epsilon$
for any $\epsilon > 0$. This agrees with the empirical evidence of
Section \ref{jack}.

We now suppose that $k \geq 7$. Let $\sigma_1 = 1, \dots, \sigma_{(k-1)/2}$
be elements of the Galois group of $\bbq(\xi_k) / \bbq$ which form a set
of representatives for the Galois group of $\bbq(\xi_k)^+ / \bbq$.
As we saw in the two previous sections,
\begin{eqnarray*}
L_E(1, \chi) = 0
&\iff& \phi(n_E(\chi)) = \left( n_E(\chi)^{\sigma_1}, \dots, n_E(\chi)^{\sigma_{(k-1)/2}} \right) \in R.
\end{eqnarray*}
As $R \subseteq R'$, and using Lemma \ref{fork7},
the probability that
$L_E(1, \chi)$ is zero is bounded by
\begin{eqnarray*}
%%\label{septicbound}
\mbox{Prob} \left( \vert L_E(1, \chi^{\sigma_1}) \vert <
\frac{c_k}{\sqrt{m}} \right), \dots, \mbox{Prob} \left( \vert L_E(1,
\chi^{\sigma_{(k-1)/2}}) \vert < \frac{c_k}{\sqrt{m}} \right) .
\end{eqnarray*}
Assuming that $|L_E(1, \chi^{\sigma_i})$  are  independent identically distributed
random variables for $1 \leq i \leq (k-1)/2$, and using (\ref{probability}), we get that the probability that
$L_E(1, \chi)$ is zero is
$$
\frac{\log^{(k-1)/8}{m}}{m^{(k-1)/4}}.
$$
neglecting all constants which are not significant here. Summing
the probabilities, this gives for $k \geq 7$
\begin{eqnarray}
\sum_{\chi \in S_3(X)}  \frac{\log^{(k-1)/8}{m}}{m^{(k-1)/4}} &=& O(1).
\end{eqnarray}
From the random matrix model, we then conjecture that the number of twists of order $k \geq 7$
such that $L_E(1, \chi)$ vanishes is bounded. This also agrees with the empirical evidence of Section \ref{jack}.

\section{Numerical Evidence}
\label{jack}

The following table shows the observed number of vanishing twists
$L_E(1, \chi)$ for characters of orders three, five and seven, and
for the first three elliptic curves in the Cremona catalogue
\cite{cremona}. For each elliptic curve $E$, the characters with
conductor prime to $N_E$ and less than two million were
considered. Any two characters of conductor $m$ and order $k$
generating the same cyclic subgroup of the character group are
conjugate, and hence the special values $L_E(1, \chi)$ vanish
simultaneously by Lemma \ref{discrete}. The number in the table
records one of each class of conjugate characters for which the
special value vanishes, which is $1/(k-1)$ of the number of
characters with vanishing special value. The twists of order
eleven in the same range for the curve $E14$ were also computed,
and no vanishing were found.

\bigskip

\centerline{
\begin{tabular}{|c|c|c|c|}
\hline Curve & Cubic & Quintic & Septic  \\
 & vanishing & vanishing & vanishing  \\
\hline
\hline E11   & 1152   &  15   & 2      \\
\hline E14   & 4347   & 10  & 0       \\
\hline E15   &  2050      & 11    & 0    \\
\hline
\end{tabular}}

\bigskip

The results for cubic twists have been analyzed in \cite{DFK} and support Conjecture \ref{mainconjecture}.
The results for quintic and septic twists are too sparse to either support or refute Conjecture \ref{mainconjecture},
but they nevertheless illustrate the extreme scarcity of vanishing in higher order twists which is
predicted by the conjecture.

\bigskip

\footnotesize{ {\bf Acknowledgments} This paper was first presented
at the ``Ranks of Elliptic Curves and Random Matrix Theory'' workshop held at the Isaac Newton
Institute in February 2004. The first and second authors would like to thank the
organizers of the workshop and the Isaac Newton Institute for their hospitality and
financial support. The first author would also
like to thank B. Birch, B. Conrey and C. Hughes for helpful discussions.}

\end{document}